\documentclass[11pt]{article}
\usepackage{amsthm,epsfig,a4}
\usepackage[ansinew]{inputenc}
\usepackage{amsmath,amssymb}
\usepackage{amsfonts}
\usepackage{eufrak}
\def\N{\mathbb{N}}
\def\R{\mathbb{R}}

\def\dr{\mathrm{d}}
\def\q{\nolinebreak \hfill $\Box$} %q.e.d.%
\def\:{\colon\thinspace}
\theoremstyle{plain}
\newtheorem{Lemma}{Lemma}[section]
\newtheorem{Prop}[Lemma]{Proposition}
\newtheorem{Thm}[Lemma]{Theorem}
\newtheorem{Cor}[Lemma]{Corollary}
\theoremstyle{definition}

\newtheorem*{Rem}{Remark}
\newtheorem*{Q}{Question}
\newtheorem*{Ac}{Acknowledgement}
\allowdisplaybreaks
\begin{document}
\title{The solution on the geography-problem of non-formal compact (almost) contact manifolds}
\author{\textsc{Christoph Bock}
%\\ \small{e-mail: bock@mi.uni-erlangen.de}
}
\date{}
\maketitle

%\begin{center}
%\textbf{Pour Iris}
%\end{center}

\begin{abstract}
Let $(m,b)$ be a pair of natural numbers. 
For $m$ odd with $m \ge 7$ (resp.\ $m \ge 5$) and $b=1$ (resp.\ $b=0$) we show that there is a non-formal compact (almost) contact $m$-manifold with first Betti number $b_1 = b$.
Moreover, in the case $b = 0$ with $m \ge 7$, the manifold even is simply-connected.
\end{abstract}
\section{Introduction}
Throughout this note a manifold is assumed to be smooth, connected and to have no boundary.
We assume that the reader is familiar with the notions of (minimal) differential graded algebras and formality as described in \cite[Section 2]{ipse1}.
''The`` \emph{minimal model of a manifold} $M$ is ''the`` minimal model for the de Rahm complex $(\Omega(M),\dr)$ of differential forms on $M$.
If the former is formal, then $M$ is called \emph{formal}.
Fern\'andez and Mu\~noz solved in \cite{FMGeo} the geography-problem of non-formal oriented compact manifolds and obtained the following theorem:

\begin{Thm}[{\cite[Theorem 1.1]{FMGeo}}] $\,$ \label{nicht sympl}
Given $m \in \mathbb{N}_+$ and $b \in \mathbb{N}$, there are oriented compact $m$-dimensional manifolds with $b_1 = b$ which are non-formal if and only if one of the following conditions holds:
\begin{itemize}
\item[(i)] $m \ge 3 \mbox{ and } b \ge 2$,

\item[(ii)] $m \ge 5 \mbox{ and } b = 1$,

\item[(iii)] $m \ge 7 \mbox{ and } b = 0$. \q
\end{itemize}
\end{Thm}

Later, we proved the analogous result for compact symplectic manifolds \cite[Theorem 1.4]{ipse1} -- of course, $m$ has to be even --, and that the next theorem holds: 

\begin{Thm}[{\cite[Theorem 1.5]{ipse1}}] $\,$ \label{Main}
For each pair $(m,b) \in \mathbb{N}_+ \times \mathbb{N}$ with $m$ odd and $b \ge 2$ exists a non-formal compact contact $m$-manifold with $b_1 = b$.
\end{Thm}

As remarked by the referee of an earlier version of this note, \cite[Theorem 3.2]{BFMT} implies the following theorem.

\begin{Thm} $\,$ \label{NewMain2}
Let $m \in \N$ with $m$ odd and $m \ge 7$.
Then there is a simply-connected non-formal compact contact $m$-manifold. \q
\end{Thm}

We shall see:

\begin{Thm} $\,$ \label{NewMain3}
There is a compact contact $5$-manifold with $b_1 = 1$ that is non-formal.
\end{Thm}

\begin{Thm} $\,$ \label{NewMain1}
Let $m \in \N$ with $m$ odd and $m \ge 7$.
Then there is a non-formal compact contact $m$-manifold with $b_1 = 1$.
\end{Thm}

\begin{Rem}
This remark is due to the referee of an earlier version of this note.
While compact Kähler manifolds are formal, see \cite[p.\ 270]{DGMS}, it happens that formality is not an obstruction to the existence of Sasakian structures on an odd-dimensional compact manifold.
This was shown in \cite[Theorem 4.9]{BFMT}. 
Thus, formality does not allow to distinguish contact manifolds which admit Sasakian structures from those which do not.
\end{Rem}

Let $n \in \N$.
Recall that a \emph{contact manifold} is a pair $(M, \xi)$, where $M$ is a $(2n+1)$-manifold and $\xi := \ker \alpha \subset TM$ a (co-orientated) hyperplane field, where $\alpha \in \Omega^1(M)$ is a differential $1$-form with $\alpha_p \wedge {(\dr \alpha)^n}_p \ne 0$ for all $p \in M$.
%A \emph{contact manifold} is a pair $(M, \xi := \ker \alpha)$, where $M$ is a $(2n+1)$-manifold and $\alpha \in \Omega^1(M)$ a differential $1$-form with $\alpha_p \wedge {(\dr \alpha)^n}_p \ne 0$ for all $p \in M$.%\footnote{Lie called the geometry of contact manifolds ,,Geometrie der Be\-rühr\-ungs\-trans\-for\-ma\-tio\-nen``.
%Meanwhile, in German, it is usual to speak of ,,Kontaktgeometrie``.
%As taught me by Geiges, this is the wrong translation of the English notion ,,contact geometry``.}
This implies that the structure group of the tangent bundle of $M$ reduces to $\{1\} \times {\mathrm U}(n)$, see \cite[p.\ 68]{GeB2}.
Therefore, it makes sence to call a $(2n+1)$-manifold \emph{almost contact}, if the structure group of the tangent bundle reduces to $\{1\} \times {\mathrm U}(n)$.
There is the following consequence of \cite[Corollary 1.3]{BER}:

\begin{Prop} $\,$ \label{almost contact contact}
Let $M$ be an almost contact manifold.
Then $M$ admits the structure of a contact manifold. \q
\end{Prop}

Therefore, the work problem is reduced to studying the geography of compact almost contact manifolds which are non-formal.

\section{Some five-dimensional contact solvmanifolds and the proof of Theorem \ref{NewMain3}}

Recall that a \emph{lattice in a Lie group} is a discrete co-compact subgroup.
Proposition \ref{almost contact contact} implies the following result:

\begin{Prop} $\,$ \label{MainP}
Let $H$ be a connected Lie group that possesses the structure of an almost complex manifold.
Then each quotient of $\R \ltimes H$ by a lattice admits a contact structure. \q
\end{Prop}

A \emph{solvmanifold} is a homogeneous space $\Gamma \backslash G$, where $G$ is a connected and simply-connected solvable Lie group and $\Gamma$ a lattice in $G$.

We are going to denote the connected and simply-connected Lie Groups as in \cite{ipse2}, i.e.\ $G_{i.j}$ is the connected and simply-connected Lie group with Lie algebra $\mathfrak{g}_{i.j}$, see \cite[Appendix A]{ipse2}.
(If any, we chose the same superscripts for the Lie groups as for their Lie algebras). 

\begin{Thm} $\,$ \label{Main5}
Let $G \notin \{G_{5.2}, \R \times G_{4.5}^{-p,-p-1}\}$, where $-\frac{1}{2} \le p < 0$, be a five-dimensional connected and simply-connected solvable Lie Group and $\Gamma$ be a lattice in $G$.
Then holds
\begin{itemize}
\item[(i)] $G \ne G_{5.2}$ is nilpotent or indecomposable non-nilpotent, or
\item[(ii)] $G = \mathbb{R} \times H$, where $H$ is a four-dimensional connected and simply-connected almost complex Lie group, or
\item[(iii)] $G = \mathbb{R} \times H$, where $H$ is a four-dimensional connected and simply-connected symplectic\footnote{A \emph{symplectic manifold} is a pair $(M, \omega)$, where $M$ is a manifold and $\omega \in \Omega^2(M)$ such taht ${{\omega}^2}_p \ne 0$ for all $p \in M$.
In this case, $\omega$ is called \emph{symplectic form on $M$.}} Lie group,
\end{itemize}
and the solvmanifold $\Gamma \backslash G$ is contact.
\end{Thm}
\pagebreak

\textit{Proof.} ad (i): By \cite[Section 7]{ipse2}, it is enough to show that the theorem holds for an almost abelian\footnote{A group is called \emph{almost abelian} if it equals $\R \ltimes \R^n$, where $n \in \N$.} connected and simply-connected Lie group $G$.
(In loc.\ cit.\ we proved that non-contact five-dimensional solvmanifolds $\Gamma \backslash G$, where $G, \Gamma$ as in (i), may only be quotients of $G_{5,2}$ or an almost abelian connected and simply-connected Lie group by a lattice.)
Proposition \ref{MainP} yields (i). 

ad (ii): \cite[Theorems 7.1.1 and 6.6.2]{ipse2} as well as Proposition \ref{MainP} imply (ii).

ad (iii): This follows from \cite[Proposition 2.4.5]{GeB2} and (ii). \q

\begin{Q} $\,$
By \cite[Section 7]{ipse2}, we considered all five-dimensional connected and simply-connected Lie groups that may possess a lattice, (such that the quotient by such a lattice could be contact) -- except for $G_{5.2}$ and $\R \times G_{4.5}^{-p,-p-1}, -\frac{1}{2} \le p < 0$.

Let $G \in \{G_{5.2}, \R \times G_{4.5}^{-p,-p-1}\}$, $-\frac{1}{2} \le p < 0$, be a five-dimensional connected and simply-connected solvable Lie group and $\Gamma$ be a lattice in $G$.
Does $\Gamma \backslash G$ admit a contact structure?
\end{Q}

\begin{Rem} $\,$
Given a solvmanifold $\Gamma \backslash G$.
Then $\Gamma$ does not determine $G$.
Note, the $3$-torus $\mathbb{Z}^3 \backslash \mathbb{R}^3$ is diffeomorphic to $\mathbb{Z}^3 \backslash G_{3.5}^0$, where $G_{3.5}^0$ equals $\mathbb{R}^3$ as a manifold and whose Lie group structure is given by
$$ (s,a,b) \cdot (t,x,y) = ( s + t, \cos(2 \pi t) \, a - \sin(2 \pi t) \, b + x, \sin(2 \pi t) \, a + \cos(2 \pi t) \, b +  y). $$ 
\end{Rem}

\textit{Proof of Theorem \ref{NewMain3}.} The non-formal $5$-manifold constructed in \cite[Proposition 7.2.9]{ipse2} is a solvmanifold with $b_1 = 1$ as a quotient of the almost abelian Lie group $G_{5.15}^{-1}$ by a lattice.
Therefore, this solvmanifold is contact by the last theorem. \q

\section{Massey products}

In order to prove non-formality, the concept of Massey products and $a$-Massey products plays an important role.
The latter were developed by Fern\'andez and Mu\~noz in \cite{FMacht}.

Let $(A,d)$ be a differential graded algebra.
\begin{itemize}
\item[(i)] Let $a_i \in H^{p_i}(A,d)$, $p_i \in \N_+, i \in \{1,2,3\}$, satisfying $a_j \cdot a_{j+1} = 0$ for $j \in \{1,2\}$.
Take elements $\alpha_i$ of $A$ with $a_i = [\alpha_i]$ and write $\alpha_j \cdot \alpha_{j+1} = d \xi_{j,j+1}$ for $j \in \{1,2\}$.
The \emph{(triple-)Massey product}\label{DMP} $\langle a_1,a_2,a_3 \rangle$ of the classes $a_i$ is defined as 
$$ [\alpha_1 \cdot \xi_{2,3} + (-1)^{p_1 +1} \xi_{1,2} \cdot \alpha_3] \in \frac{H^{p_1 + p_2 + p_3 - 1}(A,d)}{a_1 \cdot H^{p_2 + p_3 -1}(A,d) + H^{p_1 + p_2 -1}(A,d) \cdot a_3 }. $$

\item[(ii)] Let $a,b_1,b_2,b_3 \in H^2(A,d)$ satisfying $a \cdot b_i = 0$ for $i \in \{1,2,3\}$.
Take choices of representatives $a = [\alpha], b_i = [\beta_i]$ and $\alpha \cdot \beta_i = d \xi_i$ for $i \in \{1,2,3\}$. 
Then the \emph{$a$-Massey product}\label{$a$-Massey} $\langle a;b_1,b_2,b_3 \rangle$ is defined as $ [ \xi_1 \cdot \xi_2 \cdot \beta_3 + \xi_2 \cdot \xi_3 \cdot \beta_1 + \xi_3 \cdot \xi_1 \cdot  \beta_2]$ in
$$ \frac{H^8(A,d)}{\langle b_1, a, b_2 \rangle \cdot H^3(A,d) + \langle b_1, a ,b_3 \rangle \cdot H^3(A,d) + \langle b_2, a, b_3 \rangle \cdot H^3(A,d)}. $$
\end{itemize}

The next two lemmata show the relation between formality and Massey products.

\begin{Lemma}[{\cite[p.\ 260]{DGMS}}, {\cite[Theorem 1.6.5]{TO}}] $\,$
For any formal minimal differential graded algebra all Massey products vanish. \q
\end{Lemma}

\begin{Lemma}[{\cite[Proposition 3.2]{FMacht}}] $\,$
If a minimal differential graded algebra is formal, then every $a$-Massey product vanishes. \q
\end{Lemma}

\begin{Cor} $\,$
If the de Rahm complex $(\Omega(M),\dr)$ of a manifold $M$ possesses a non-vanishing Massey or $a$-Massey product, then $M$ is not formal. \q
\end{Cor}

We will need the following lemma, too.

\begin{Lemma}[{\cite[Lemma 2.11]{FMDon}}] \label{Prod}
The product of two manifolds is formal if and only if both factors are formal. \q
\end{Lemma}
\section{Proof of Theorem \ref{NewMain1}}
In \cite[Theorem 8.3.2]{ipse2}, we considered the completely solvable Lie Group $G := G_{6.15}^{-1}$ and constructed a lattice $\Gamma$ in $G$.
The space $\mathfrak{g}^*$ of left-invariant differential $1$-forms on $G$ possesses a basis $\{x_1, \ldots, x_6\}$ such that
$$\dr x_1 = - x_{23}, \, \dr x_2 = - x_{26}, \, \dr x_3 = x_{36}, \, \dr x_4 = -x_{26} - x_{46} , \, \dr x_5 = -x_{36} + x_{56}, \, \dr x_6 = 0$$
and therefore (for the non-exact $x_{ij}$)
$$ \begin{array}{l@{~~~~}l@{~~~~}l}
\dr x_{12} = x_{126}, & \dr x_{13} = -x_{136}, & \dr x_{14} = x_{126} + x_{146} - x_{234}, \\
\dr x_{15} = x_{136} - x_{156} - x_{235}, & \dr x_{16} = -x_{236}, & \dr x_{24} = 2 x_{246}, \\
\dr x_{25} = x_{236}, & \dr x_{34} = - x_{236}, & \dr x_{35} = -2 x_{356}, \\
\dr x_{45} = x_{256} - x_{346}, &&
\end{array} $$
where $x_{ij} := x_i \wedge x_j$ as well as $x_{ijk} := x_i \wedge x_j \wedge x_k$, $i,j,k \in \{1, \ldots, 6\}$.%
\footnote{There is a misprint in the proof of \cite[Theorem 8.3.2 (ii)]{ipse2}.
There, $\delta x_{16}$ has to equal $-x_{236}$ instead of $x_{236}$.} %
This yields
$$ \dr x_{123} = \dr x_{256} = \dr x_{346} = \dr x_{456} = 0, $$
$$ \begin{array}{l@{~~~~}l@{~~~~}l}
\dr x_{124} = -2 x_{1246}, & \dr x_{125} = -x_{1236}, &\dr x_{134} = x_{1236}, \\
\dr x_{135} = 2 x_{1356}, & \dr x_{145} = -x_{1256} + x_{1346} - x_{2346}, & \dr x_{146} = - x_{2346}, \\
\dr x_{156} = - x_{2346}, & \dr x_{234} = - x_{2346}, & \dr x_{235} = x_{2356}, \\
\dr x_{245} = x_{2346} - x_{2356}, & \dr x_{345} = x_{2356} + x_{3456}, &
\end{array} $$
where $x_{ijkl} := x_i \wedge x_j \wedge x_k \wedge x_l$, $i,j,k,l \in \{1, \ldots, 6\}$, and
$$ \begin{array}{l@{~~~~}l@{~~~~}l}
\dr x_{1234} = x_{12346}, & \dr x_{1235} = -2 x_{12356}, & \dr x_{1245} = -x_{12346} + x_{12456}, \\ 
\dr x_{1256} = 0, & \dr x_{1345} = x_{12356} - x_{13456}, & \dr x_{1346} = 0, \\
\dr x_{1456} = - x_{23456}, & \dr x_{2345} = 0, & \dr x_{12345} = 0,
\end{array} $$
where $x_{ijklm} := x_i \wedge x_j \wedge x_k \wedge x_l \wedge x_m$, $i,j,k,l,m \in \{1, \ldots, 6\}$.

By completely solvability, \cite[Theorem 3.2.10]{TO} and Hattori's Theorem \cite[p.\ 77]{TO}\footnote{These two theorems were quoted in \cite[Theorem 3.10 (i), (ii)]{ipse2}.}, the cohomology groups of the corresponding solvmanifold $M := \Gamma \backslash G$ are given by
\begin{gather}
H^1(M,\R) \cong \langle [x_6] \rangle_{\R}, ~~ H^2(M,\R) \cong \langle [x_{16} + x_{25}], [x_{16} - x_{34}] \rangle_{\R}, \label{M.1.1} \\
H^3(M,\R) \cong \langle [x_{123}], [x_{125} + x_{134}], \underbrace{[x_{256}]}_{= - [x_{346}]}, [x_{456}] \rangle_{\R}, \label{M.1.2} \\
H^4(M,\R) \cong \langle \underbrace{[x_{1256}]}_{= [x_{1346}]}, [x_{2345}] \rangle_{\R}, \label{M.1.3} \\
H^5(M,\R) \cong \langle [x_{12345}] \rangle_{\R} \nonumber.
\end{gather} 

$\tilde{\omega} := 2 x_{16} + x_{25} - x_{34} = (x_{16} + x_{25}) + (x_{16} - x_{34})$ induces a symplectic form on $M$.
Analogous to \cite[Observation 4.3]{Go}, one can see that there is a symplectic form $\omega$ on $M$ with
\begin{equation}
[\omega] = \lambda \, [x_{16} + x_{25}] + \mu \, [x_{16} - x_{34}], ~ \lambda, \mu \in \R \setminus \{0\}, \lambda \ne -\mu, \label{M.1.0}
\end{equation}
and $[\omega]$ lifts to an integral cohomlogy class -- namely $\lambda = \frac{1 + \varepsilon_1}{n_0}, \mu = \frac{1 + \varepsilon_2}{n_0}$ for certain $\varepsilon_1, \varepsilon_2  \in \R_+, n_0 \in \N_+$.
By Boothby and Wang \cite[Theorem 3]{BW}, there is a principal fibre bundle
$$ S^1 \longrightarrow E \stackrel{\pi}{\longrightarrow} M, $$
where $E$ is a compact contact manifold.
We apply the Gysin sequence
$$ \begin{array}{ccccccc}
\{0\} & \longrightarrow                           & H^1(M,\R) & \stackrel{\pi^*}{\longrightarrow} & H^1(E,\R) & \longrightarrow & H^0(M,\R) \\
      & \stackrel{[\omega] \cup}{\longrightarrow} & H^2(M,\R) & \stackrel{\pi^*}{\longrightarrow} & H^2(E,\R) & \longrightarrow & H^1(M,\R) \\
      & \stackrel{[\omega] \cup}{\longrightarrow} & H^3(M,\R) & \stackrel{\pi^*}{\longrightarrow} & H^3(E,\R) & \longrightarrow & H^2(M,\R) \\
      & \stackrel{[\omega] \cup}{\longrightarrow} & H^4(M,\R) & \longrightarrow                   & \ldots    & 
\end{array} $$
to obtain
\begin{gather}
H^1(E,\R) \cong \langle \pi^* [x_6] \rangle_{\R}, \label{M.1.4} \\
H^2(E,\R) \cong \langle \pi^* [x_{16} + x_{25}], \pi^* [x_{16} - x_{34}] \, | \, \lambda \, \pi^* [x_{16} + x_{25}] + \mu \, \pi^* [x_{16} - x_{23}] = 0 \rangle_{\R}, \label{M.1.5} \\
H^3(E,\R) \cong \langle \pi^* [x_{123}], \pi^*[x_{125} + x_{134}], \pi^* [x_{456}] \rangle_{\R}. \label{M.1.6}
\end{gather} 

{[} Obviously, (\ref{M.1.4}) holds.

Using (\ref{M.1.0}), (\ref{M.1.1}), (\ref{M.1.2}), one easily shows the injectivity of
$$ [\omega] \cup \ldots \: H^1(M,\R) \longrightarrow H^3(M,\R), $$
i.e.\ the surjectivity of
\begin{equation}
\pi^* \: H^2(M,\R) \longrightarrow H^2(E,\R). \label{M.1.7}
\end{equation}
Since the kernel of (\ref{M.1.7}) equals $\langle [\omega] \rangle_{\R}$, (\ref{M.1.5}) holds by (\ref{M.1.0}).

Analogous, one computes that
$$ [\omega] \cup \ldots \: H^2(M,\R) \longrightarrow H^4(M,\R) $$
is injective.
So,
\begin{equation}
\pi^* \: H^3(M,\R) \longrightarrow H^3(E,\R) \label{M.1.8}
\end{equation}
is surjective with kernel $\langle [x_{256}] \rangle_{\R}$, and (\ref{M.1.6}) follows. {]}

We have $x_{66} = \dr 0$, $x_6 \wedge \tilde{\omega} = x_{256} - x_{346} = \dr x_{45}$ as well as $x_6 \wedge x_{45} + 0 \wedge \tilde{\omega} = x_{456}$ and $x_{456}$ is not exact.
Since
$$ [x_6] \cup H^2(M,\R) + H^1(M,\R) \cup [\tilde{\omega}] \stackrel{(\ref{M.1.1})}{=} \langle [x_{256}], [x_{346}] \rangle_{\R} \stackrel{(\ref{M.1.2})}{=} \langle [x_{256}] \rangle_{\R}, $$
$\langle [x_6], [x_6], [\tilde{\omega}] \rangle = [x_{456}] \in \langle [x_{123}], [x_{125} + x_{134}], [x_{456}] \rangle_{\R}$ is a non-vanishing Massey pro\-duct.
(\ref{M.1.4}) -- (\ref{M.1.6}) imply that $\langle \pi^* [x_6], \pi^* [x_6], \pi^* [\tilde{\omega}] \rangle = \pi^* [x_{456}] $ also is a non-vanishing Massey product.
Note,
$$ \pi^* [x_6] \cup H^2(E,\R) + H^1(E,\R) \cup \pi^* [\tilde{\omega}] \stackrel{(\ref{M.1.4}), (\ref{M.1.5}), (\ref{M.1.2})}{=} \langle \pi^* [x_{256}] \rangle_{\R} = \{0\} $$
because the kernel of (\ref{M.1.8}) equals $\langle [x_{256}] \rangle_{\R}$.

The case $m = 7$ follows.
For $m = 2n + 1$ with $n \in \N_+, \, n \ge 4,$ consider the man\-i\-folds $M \times (S^2)^{(n-3)}$ instead of $M$.
Clearly, these manifolds are symplectic.
By Lemma \ref{Prod}, they are non-formal. \q

\section{Another proof of Theorem \ref{NewMain2} in the case $m \ge 13$}

We give an idea for a shorter proof of the theorem in the case $m \ge 13$:
Cavalcanti proved in \cite[Example 4.4]{C} the existence of a simpy-connected non-formal compact symplectic $12$-manifold $M$ that has a non-vanishing Massey product coming from three differential $1$-forms.
\cite[Proposition 6.3]{ipse1} yields the theorem for $m = 13$.
Again, by considering the product of $M$ with finitely many copies of $S^2$, one obtains the higher-dimensional examples. \q

\begin{Rem} $\,$
Fern\'andez and Mu\~noz constructed in \cite{FMacht} an $8$-dimensional non-formal compact symplectic manifold $(M,\tilde{\omega})$ with
\begin{equation*} \label{M.2.1}
\begin{array}{l@{~~~}l@{~~~}l}
b_0(M) = b_8(M) = 1, &  b_1(M) = b_7(M) = 0, & b_2(M) = b_6(M) = 256, \\
b_3(M) = b_5(M) = 0, & b_4(M) = 269. & \\
\end{array}
\end{equation*}
There is an $a$-Massey product $\langle [\alpha]; [\beta_1], [\beta_2], [\beta_3] \rangle$ for certain closed $2$-forms $\alpha, \beta_i$, where $i \in \{1,2,3\}$, on $M$:
One has $\langle [\alpha]; [\beta_1], [\beta_2], [\beta_3] \rangle = \lambda \,[ \tilde{\omega}^4 ]$ for $\lambda \ne 0$.
Clearly, $\lambda \, \tilde{\omega}^4$ is not exact, and since $b_3(M) = 0$, it follows that this $a$-Massey product does not vanish.
Again, by \cite[Observation 4.3]{Go}, there is a symplectic form $\omega$ on $M$ whose cohomology class lifts to an integral cohomolgy class, and we have the Boothby-Wang fibration 
$$ S^1 \longrightarrow E \stackrel{\pi}{\longrightarrow} M, $$
where $E$ is a compact contact manifold with $\dim E = 9$. 
The Gysin sequence yields $H^1(E,\R) = \{0\}$, i.e. $H^8(E,\R) = \{0\}$.
Therefore, $\langle \pi^* [\alpha]; \pi^* [\beta_1], \pi^* [\beta_2], \pi^* [\beta_3] \rangle$ vanishes.
But we do not know whether $E$ is non-formal.
\end{Rem}

\begin{Ac}
The author wishes to thank Marisa Fern\'andez and Hansjörg Geiges for various conversations with both of them.
\end{Ac}

\noindent
\textsc{Christoph Bock\\ Department Mathematik\\ Universität Erlangen-N\"urn\-berg\\ Cauerstraße 11\\ 91058 Erlangen\\ Germany}

\noindent
\textit{e-mail:} \verb"bock@mi.uni-erlangen.de"
\end{document}